\documentclass[12pt]{article}
\usepackage{graphicx, amsmath}
\usepackage{amsfonts}
\usepackage{verbatim}
\usepackage{latexsym}
\makeatletter
\@addtoreset{equation}{section}

\begin{document}

\newcommand{\E}{\mathbb{E}}
\newcommand{\PP}{\mathbb{P}}
\newcommand{\RR}{\mathbb{R}}

\newtheorem{theorem}{Theorem}[section]
\newtheorem{lemma}{Lemma}[section]
\newtheorem{coro}{Corollary}[section]
\newtheorem{defn}{Definition}[section]
\newtheorem{assp}{Assumption}[section]
\newtheorem{expl}{Example}[section]
\newtheorem{prop}{Proposition}[section]
\newtheorem{remark}{Remark}[section]

\newcommand\tq{{\scriptstyle{3\over 4 }\scriptstyle}}
\newcommand\qua{{\scriptstyle{1\over 4 }\scriptstyle}}
\newcommand\hf{{\textstyle{1\over 2 }\displaystyle}}
\newcommand\hhf{{\scriptstyle{1\over 2 }\scriptstyle}}

\newcommand{\eproof}{\indent\vrule height6pt width4pt depth1pt\hfil\par\medbreak}

\def\a{\alpha} \def\g{\gamma}  
\def\e{\varepsilon} \def\z{\zeta} \def\y{\eta} \def\o{\theta}
\def\vo{\vartheta} \def\k{\kappa} \def\l{\lambda} \def\m{\mu} \def\n{\nu}
\def\x{\xi}  \def\r{\rho} \def\s{\sigma}   
\def\p{\phi} \def\f{\varphi}   \def\w{\omega}  
\def\q{\surd} \def\i{\bot} \def\h{\forall} \def\j{\emptyset}

\def\be{\beta} \def\de{\delta} \def\up{\upsilon} \def\eq{\equiv} 
\def\ve{\vee} \def\we{\wedge} 

\def\F{{\cal F}}
\def\T{\tau} \def\G{\Gamma}  \def\D{\Delta} \def\O{\Theta} \def\L{\Lambda}
\def\X{\Xi} \def\Si{\Sigma} \def\W{\Omega} 
\def\M{\partial} \def\N{\nabla} \def\Ex{\exists} \def\K{\times}
\def\V{\bigvee} \def\U{\bigwedge}

\def\1{\oslash} \def\2{\oplus} \def\3{\otimes} \def\4{\ominus} 
\def\5{\circ} \def\6{\odot} \def\7{\backslash} \def\8{\infty} 
\def\9{\bigcap} \def\0{\bigcup} \def\+{\pm} \def\-{\mp}
\def\la{\langle} \def\ra{\rangle}

\def\tl{\tilde} 
\def\trace{\hbox{\rm trace}}
\def\diag{\hbox{\rm diag}}
\def\for{\quad\hbox{for }}
\def\refer{\hangindent=0.3in\hangafter=1}

\newcommand\wD{\widehat{\D}}
\newcommand{\ka}{\kappa_{10}}

\title{
\bf {Neutral Stochastic Differential  Delay Equations with Locally Monotone Coefficients}
}
\author
{ {\bf Yanting Ji}$^{\tt a}$ \, {\bf Qingshuo Song}$^{\tt b}$ \,
{\bf Chenggui Yuan}$^{\tt a}$\thanks{Contact e-mail address: C.Yuan@swansea.ac.uk}\
\\[1ex]
$^{\tt a}$Department of Mathematics\\
University of Wales Swansea, Swansea,  SA2 8PP, U.K.\\
$^{\tt b}$Department of Mathematics\\
City University of Hong Kong, Hong Kong\\}

\date{}

\maketitle

\begin{abstract}
In this paper,  we prove the existence and uniqueness of the solution  for neutral stochastic differential delay equations  with locally monotone coefficients by using numerical approximation. An example is provided to illustrate our theory.\\

 \smallskip


\end{abstract}

\section{Introduction}

The theory of stochastic functional differential equations (SFDEs) has been developed for a while, for instant \cite{MM84} provides systematic presentation for the existence and uniqueness, Markov property, the generator and the regularity of the solutions  of SFDEs. \cite{Mao} presents the estimation of the moment of the solutions, in particular, the  Razumikhin theorem was generalized from functional differential equations to SFDEs.   For the studies of long-term behaviour of SFDEs, we here only mention \cite{ Es, HMS, RR}.

On the other hand, most SFDEs can not be solved explicitly,  numerical methods become one of the most powerful tools tackling these problems in the real world practise. There is extensive literature in investigating the strong
convergence, weak convergence or sample path convergence  of
numerical schemes for SFDEs, we here highlight  \cite{Hu, hu96, jwy09,KP00,ms03}, to name a few. 

More recently, a class of stochastic equations has emerged, which depends on the past and present values but that involves derivatives with delays as well as the function itself. Such equations are called neutral stochastic functional differential equations (NSFDEs). The theory of NSFDEs has recently received a lot of attention.  For example,   the existence and uniqueness, and the stability  of the solutions of NSFDEs can be found in \cite{Mao}. For the approximation and  numerical solutions in this area the reader is refer to \cite{Bou, Zhou}.  For large
deviation of functional NSFDEs, we refer to 
\cite{BY}.

The existence and uniqueness of solutions of stochastic equations is always an important topic. It is interesting that Krylov \cite{KR} gave a theorem for the existence and uniqueness by Euler numerical approximation under local monotonicity condition, which is much weaker than global Lipschitz condition. Recently Gy\"ogy and Sabanis \cite{Sab} extended this result to stochastic differential delay equations (SDDEs). However, up to our best knowledge,  we do not know if  neutral stochastic differential delay equations  (NSDDEs) has a unique solution under a local monotonicity condition. The main aim of this paper is to  fill the gap by extending the existed methods to establish the existence and uniqueness theorem of NSDDEs.

Throughout this paper, let $(\Omega,\mathcal{F},\mathbb{P})$ be a complete probability space with a filtration $\{\mathcal{F}_t\}_{t\geq 0}$ satisfying the usual condition (i.e. it is right continuous and $\mathcal{F}_{0}$ contains all $P$-null sets). Let $|\cdot|$ denote the Euclidean norm and $\|\cdot\|$ the matrix trace norm and $\langle\cdot,\cdot\rangle$ denotes the Euclidean inner product on $\mathbb{R}^{n}$. Let $b>a$ be two real constants and $C([a, b];\mathbb{R}^{n})$ the space of all continuous function from $[a, b]$ to $\mathbb{R}^{n}$ with the norm $\|\phi\|_{(a,b)}=\sup_{a\leq\theta\leq b}|\phi(\theta)|$. Denote by $C^b_{\mathcal{F}_{0}}([-\tau,0];\mathbb{R}^{n})$ the family of all bounded, $\mathcal{F}_{0}$-measurable, $C([-\tau,0];\mathbb{R}^{n})$-valued random variables. Denote $\mathcal{L}^{p}([a,b];\mathbb{R}^{n})$ the family of $\mathbb{R}^{n}$-valued $\mathcal{F}_{t}-$adapted process$\{h(t)\}_{a\leq t\leq b}$ such that$\int_{a}^{b}|h(t)|^pdt<\infty\quad a.s.$ Let $B(t)$ be a standard $m$-dimensional Brownian motion. Denote $C$ a generic positive constant, whose value may change from line to line.

Let $f(x, t, \w)$ and $g(x, t, \w)$ be  given as follows:
\begin{align*}
f:  \RR^n\K [0, \8)\K \W\to \RR^{n},\quad \quad
g: \RR^n\K [0, \8)\K \W\to \RR^{n\K m}
\end{align*}
such that both are continuous in $x\in \RR^n$ for each fixed $t\in [0, \8),$  and progressively measurable. In particular, for every $x\in \RR^n, t\in [0, \8)$ both are $\F_t-$measurable.
 We also assume the following conditions:
\begin{enumerate}
\item[(i)] For every $R>0, T>0,$ 
$$
\int_0^T\sup_{|x|\le R}\{ |f(x, t)|+\|g(x, t)\|^2\}dt <\8.
$$
\item[(ii)] For every $R>0, t>0, |x|\le R, |y|\le R,$  there exist  $M_R(t), M(t)$ such that $M_R(t), M(t)\in \mathcal{L}^{1}([0 , T];\mathbb{R}), \forall T>0$ and
\begin{align}
&2\la x-y, f(x, t)-f(y, t)\ra +\|g(x, t)-g(y, t)\|^2\le M_R(t)|x-y|^2,\\
& 2\la x, f(x, t)\ra+\|g(x, t)\|^2\le M(t)(1+|x|^2).
\end{align}

\end{enumerate}
Consider 
\begin{equation}\label{Kr1}
dZ(t)=f(Z(t), t)dt+g(Z(t), t)dB(t), \quad t \in [0, T],
\end{equation}
with an initial value $Z(0)$ which is $\F_{0}-$measurable. The  the following result was proved by Krylov \cite{KR}, also see Pr\'ev\^ot and R\"ockner \cite{Roc}.
\begin{theorem}\label{Kr}
Under the assumptions (i) and (ii), the equation \eqref{Kr1} has a unique solution.
\end{theorem}
Recently, Gy\"ongy and Sabanis \cite{Sab} extended this result to SDDEs. 
In this paper,  we shall generalize the existing results and extend to the NSDDEs case.
Consider an $n$-dimensional NSDDE of the following form
\begin{equation}\label{NSDDE}
d[X(t)-D(X(t-\tau))]  = b(X(t),X(t-\tau),t)dt+\sigma(X(t),X(t-\tau),t)dB(t)
\end{equation}
on $t \geq 0$. We  assume
\begin{equation}
\begin{split}
& D:\mathbb{R}^{n} \to \mathbb{R}^{n},\quad
 b:\mathbb{R}^{n}\times \mathbb{R}^{n}\K [0,T] \to \mathbb{R}^{n},\\
&\sigma:\mathbb{R}^{n}\times \mathbb{R}^{n}\K [0,T] \to \mathbb{R}^{n\times m}.
\end{split}
\end{equation}
We  also assume that $D$, $b$ and $\sigma$ are Borel-measurable and the initial data is given by:
\begin{equation}\label{initial}
X(0) = \xi \in C^b_{\mathcal{F}_{0}}([-\tau,0];\mathbb{R}^{n}).
\end{equation}
Throughout the paper, for  $T>\tau>0$ we assume that $T/\tau$ is a rational number.
%Now, Fix a $T \geq 0$, without loss of generality, we T/\tau assume that $T = k\cdot\tau,$ where $k \in \mathbb{N}.$\\

Firstly, we need to impose the standing integrability hypothesis for this paper:
\begin{description}
\item[(H)] For every $R > 0, T>0$  
\begin{equation}\label{Hypothesis}
\int_0^T\sup_{|x|\le R, |y|\le  R}\{|b(x,y,t)|+ \|\sigma(x,y,t)\|^2\}dt <\8 \mbox{ on } \Omega.
\end{equation}
\end{description}
We now give the definition of the solution to the equation $(\ref{NSDDE})$  with initial data $(\ref{initial}).$
\begin{defn}
An $\mathbb{R}^{n}$-valued stochastic process $X(t)$ on $[-\tau,T]$ is called a solution to equation $(\ref{NSDDE})$ with initial data $(\ref{initial})$ if it has following properties:
\begin{enumerate}
\item[\rm{(i)}] It is continuous and $\{X(t)_{0\leq t \leq T}\}$ is $\mathcal{F}_{t}$-adapted.
\item[\rm{(ii)}] $\{b(X(t),X(t-\tau),t)\} \in \mathcal{L}^{1}([0,T];\mathbb{R}^{n})$ and $\{\sigma(X(t),X(t-\tau),t)\} \in \mathcal{L}^{2}([0,T];\mathbb{R}^{n\times m})$
\item[\rm{(iii)}] $X(0) = \xi(0)$ and 
\begin{equation}
\begin{split}
X(t) - D(X(t-\tau))& =  X(0) - D(\xi(-\tau))+\int_{0}^t b(X(s),X(s-\tau),s)ds\\
& +\int_{0}^t\sigma(X(s),X(s-\tau),s)dB(s)
\end{split}
\end{equation}
hold with probability one, for all $0 \leq t \leq T.$\\
A solution $X(t)$ is said to be unique if any other solution $\bar{X}(t)$ is indistinguishable from it, that is 
\begin{equation*}
P\{X(t) =\bar{X}(t) \quad \mbox{for all} \quad -\tau \leq t \leq T \} = 1.
\end{equation*}
\end{enumerate}
\end{defn}
The following theorem is our main result.
\begin{theorem}\label{EUtheorem}
Assume   $D$, $b$ and $\sigma$ 
satisfy the following assumptions for all $T,R\in [0,\infty)$:
\begin{description}
\item[(C1)] The functions $b(x,y,t), \sigma(x, y, t)$ are continuous in both $x$ and $y$ for all $t\in [0,T].$ 
\item[(C2)]  There exist two $\mathbb{R}_{+}-$ valued 
functions $K_{1}(t),$ $\tilde{K}_{1}(t)$ and a positive constant $C_1(\tau)$ such that for all $t\in [0,T],$
 $K_{1}(t) \leq  C_{1}(\tau)K_{1}(t-\tau),$ $K_{1}(t)\geq \tilde{K}_{1}(t)$ and 
\begin{equation}\label{WeakCorecivity}
2\langle x-D(y),b(x,y,t)\rangle + \|\sigma(x,y,t)\|^{2}\leq K_{1}(t)(1+|x|^{2})+\tilde{K}_{1}(t-\tau)(1+|y|^2),
\end{equation}
for $\forall$ $x,y \in \mathbb{R}^{n},$ $t\in [0,T].$\\
\item[(C3)]There exist two $\mathbb{R}_{+}-$ valued functions $K_{R}(t)$, $\tilde{K}_{R}(t)$  and a positive constant $C_R(\tau)$such that for all  $t\in [0,T],$ $K_{R}(t)\leq C_{R}(\tau)K_{R}(t-\tau),$ $K_{R}(t)\geq \tilde{K}_{R}(t)$ and
\begin{equation}\label{LocalLipschitzII}
\begin{split}
&2\langle x-D(y)-\bar{x}+D(\bar{y}), b(x,y,t)-b(\bar{x},\bar{y},t)\rangle+\|\sigma(x,y,t)-\sigma(\bar{x},\bar{y},t)\|^2 \\
&\leq K_{R}(t) |x-\bar{x}|^{2}+\tilde{K}_{R}(t-\tau)|y-\bar{y}|^2 
\end{split}
\end{equation}
for all $|x|\vee|y|\vee|\bar{x}|\vee|\bar{y}| \leq R,$ $t\in [0,T].$\\
\item[(C4)] Assume $D(0) = 0$ and that there is a constant $\kappa \in (0,1)$ such that
\begin{equation*}
|D(x)-D(y)| \leq \kappa |x-y|,
\end{equation*}
holds for all $x,y\in \mathbb{R}^{n}.$
\end{description}
Moreover, we assume  $C_{1}(\tau)\vee C_{R}(\tau)\leq \frac{1}{\kappa}$ and  $K_{1}(t), \tilde{K}_{1}(t), K_{R}(t), \tilde{K}_{R}(t) \in \mathcal{L}^{1}([-\tau, T];\mathbb{R}_+).$ 
Then there exists a unique process $\{X(t)\}_{t\in [0,T]}$ that satisfies equation $(\ref{NSDDE})$ with the initial data $(\ref{initial}).$  Moreover, the mean square of the solution is finite.
\end{theorem}

\begin{remark}
If $D\equiv 0, $ then the  equation $(\ref{NSDDE})$ becomes a SDDE, which has been investigated in \cite{Sab}.  However our conditions are weaker than those in  \cite{Sab}, since the conditions in present paper include the delay components at the right hand side of \eqref{WeakCorecivity} and \eqref{LocalLipschitzII}. Moreover, If $D\equiv 0, \tau=0,$ we  take $\tilde K_1(t)=\tilde K_R(t)=0, C_1(\tau)=C_R(\tau)=1,$ then Theorem \ref{EUtheorem} becomes Theorem \ref{Kr}, which means our result is a generalization of Krylov result Theorem \ref{Kr}.
\end{remark}

The remainder of this paper is organized as follows. In Section 2 we shall give a localization lemma, which will be crucial for the proof of the
main result  Theorem \ref{EUtheorem}.  In Section 3 the proof of the main result will be demonstrated. An illustrative example will be presented in the Section 4.

\section{Localization Lemma}

In preparation for the proof of main result, Theorem \ref{EUtheorem}, we need to introduce following lemmas.
\begin{lemma}\label{lemma1}  Let $Y(t),$ $t \in [0,T],$ be a continuous, $\mathbb{R}_{+}$-valued, $\mathcal{F}_{t}$-adapted process on $(\Omega,\mathcal{F},\mathbb{P})$ and $\tau$ be a $\mathcal{F}_{t}$-stopping time, and let $\epsilon \in (0,\infty).$ Denote 
\begin{equation*}
\rho_{\epsilon} := \rho \wedge \inf \{t \geq 0 | Y(t) \geq \epsilon\}, 
\end{equation*}
then
\begin{equation*}
\mathbb{P}(\{\sup_{t \in [0,\rho]}  Y(t) \geq \epsilon\}) \leq \frac{1}{\epsilon}\E(Y(\rho_{\epsilon})).
\end{equation*}
\end{lemma}

\begin{lemma}\label{lemma2} 
Let $p >1,$ $\epsilon >0$ and $a,b \in R.$ Then
\begin{equation}
|a+b|^p \leq [1+\epsilon^{\frac{1}{p-1}}]^{p-1}\bigg(|a|^p+\frac{|b|^p}{\epsilon}\bigg).
\end{equation}
\end{lemma}

\begin{lemma}\label{LemmaKey}
Let $p>1$ and condition {\bf(C4)} hold. Then
\begin{equation*}
\sup_{0\leq s \leq t}|X(s)|^{p} \leq \frac{\kappa}{1-\kappa}||\xi||^{p}+\frac{1}{(1-\kappa)^{p}}\sup_{0\leq s \leq t}|X(s)-D(X(s-\tau))|^{p}.
\end{equation*}
\end{lemma}

The proof of lemma \ref{lemma1} can be found in \cite{Roc},
% Lemma 3.1.3 on p.44, 
the proof of Lemma \ref{lemma2} 
%Lemma 4.1 on p.211, 
and \ref{LemmaKey} can be found in \cite{Mao}. 
%Lemma 4.4 on p.212.
The following lemma is an extended version of Lemma 3.1.4 in \cite{Roc} to NSDDEs.  Since the neutral term and the delay variables are involved, the proof of following lemma is much more technical.
\begin{lemma}\label{lemma3}
Let $n \in \mathbb{N}$ and $X_{n}(t),$ $t\in [-\tau,T],$ be a continuous, $\mathbb{R}^{n}$-valued, $\mathcal{F}_{t}$-adapted process on $(\Omega,\mathcal{F},\PP)$ such that for $t\in [-\tau,0],$ $X_{n}(t) = \xi(t).$ For $t \in [0,T]$ 
\begin{equation}
\begin{split}
d[X_{n}(t)-D(X_{n}(t-\tau))] &= b(X_{n}(t)+p_{n}(t),X_{n}(t-\tau)+p_{n}(t-\tau),t)dt\\
&+\sigma(X_{n}(t)+p_{n}(t),X_{n}(t-\tau)+p_{n}(t-\tau),t)dB(t),
\end{split}
\end{equation}
for some progressively measurable process $p_{n}(t),$ and $p_{n}(t) = 0,$ for any $t\in[-\tau,0].$ For $n \in \mathbb{N}$ and $R \in [0,\infty),$ let $\tau_{n}(R)$ be $\mathcal{F}_{t}$- stopping times such that
\begin{description}
\item[(i)]
\begin{equation*}
|X_{n}(t)| + |p_{n}(t)| \leq R \quad \mbox{if} \quad t \in [0,\tau_{n}(R)] \quad \mbox{a.s.}
\end{equation*}
\item[(ii)] 
\begin{equation*}
\lim_{n \to \infty}\E\int_{0}^{T\wedge \tau_{n}(R)}|p_{n}(t)|dt = 0 \quad \mbox{for all} \quad T \in [0,\infty).
\end{equation*}
\item[(iii)] 
Assume that there exists a function $r :[0,\infty) \rightarrow [0,\infty)$ such that
$\lim_{R\rightarrow\infty}r(R)=\infty.$ Also assume that: 
\begin{equation}
\lim_{R\rightarrow\infty}\overline{\lim_{n\rightarrow\infty}}\PP\left(\left\{\tau_{n}(R)\leq T,\ \sup_{t\in[0,\tau_{n}(R)]}|X_{n}(t)-D(X_{n}(t-\tau))|<r(R)\right\}\right)=0.
\end{equation}
for all $T\in[0,\infty).$
\end{description}
Then for every $T \in [0,\infty)$ we have 
\begin{equation*}
\sup_{t\in [0,T]}|X_{n}(t)-X_{m}(t)| \xrightarrow{\mathbb{P}} 0
\end{equation*}
as $n,m \to \infty.$
\end{lemma}
{\bf Proof}:\quad We divide the proof into three steps:

{\it Step (i)} \quad
By $(\ref{Hypothesis})$ we may assume that
\begin{equation}\label{bound}
\sup_{|x|\le R, |y|\le R}|b(x,y,t)| \leq K_{R}(t),
\end{equation}
otherwise, we replace $K_{R}(t)$ by the maximum of $K_{R}(t)$ and the integrand in \eqref{Hypothesis}.
Fix $R \in [0,\infty),$ define a $\mathcal{F}_{t}$-stopping time
\begin{equation*}
\tau(R,u):= \inf\{t \geq 0 |\alpha_{R}(t) > u\}, \quad u \in [0,\infty),
\end{equation*}
where $\alpha_{R}(t): = \int_{0}^{t}K_{R}(s)ds.$
Since $K_{R}(t)\in \mathcal{L}^{1}([0 , T]$ for ny $T\ge 0$, then  $\tau(R,u)\uparrow \infty$ as $u \to \infty.$
In particular, there exist $u(R) \in [0,\infty)$ such that
\begin{equation*}
\PP(\{\tau(R,u(R)) \leq R \}) \leq \frac{1}{R}.
\end{equation*}
Now, denote $\tau(R):= \tau(R,u(R)),$ we have $\tau(R) \to \infty$ in probability as $R \to \infty$ and $\alpha_{R}(t\wedge\tau(R)) \leq u(R)$ for all $t\in [0,T]$ and $R \in [0,\infty).$
Moreover, if we replace $\tau_{n}(R)$ by $\tau_{n}(R) \wedge \tau(R)$ for $n \in \mathbb{N}$ and $R \in [0,\infty),$ we still have 
\begin{equation}
|X_{n}(t)|+|p_{n}(t)|\leq R \mbox{\ \ if\ \ }
t\in[0,\tau_{n}(R) \wedge \tau(R)]\quad \mbox{a.s.},
\end{equation}
and
\begin{equation}
\lim_{n\rightarrow\infty}\mathbb{E}\int^{T\wedge
\tau_{n}(R) \wedge \tau(R)}_{0}|p_{n}(t)|dt=0 \mbox{\ \ for\ all\ \ }
T\in[0,\infty).
\end{equation}
i.e. assumptions {\bf (i)} and {\bf (ii)} hold. 
Meanwhile, we have 
\begin{equation*}
\begin{split}
& \PP\left(\left\{\tau_{n}(R)\wedge\tau(R)\leq T,\sup_{t\in[0,\tau_{n}(R)\wedge \tau(R)]}|X_{n}(t)-D(X_{n}(t-\tau))|\leq r(R)\right\}\right)\\
&=\PP\left(\left\{\tau_{n}(R)\leq T,\sup_{t\in[0,\tau_{n}(R)]}|X_{n}(t)-D(X_{n}(t-\tau))|\leq r(R),\tau_{n}(R)\leq\tau(R)\right\}\right)\\
&+\PP\left(\left\{\tau(R)\leq T,\sup_{t\in[0,\tau(R)]}|X_{n}(t)-D(X_{n}(t-\tau))|\leq r(R),\tau_{n}(R)>\tau(R)\right\}\right)\\
&\leq \PP\left(\left\{\tau_{n}(R)\leq
T,\sup_{t\in[0,\tau_{n}(R)]}|X_{n}(t)-D(X_{n}(t-\tau))|\leq r(R),\tau_{n}(R)\leq \tau(R)\right\}\right)\\
&+\PP\left(\left\{\tau(R)\leq 
T,\tau_{n}(R)>\tau(R)\right\}\right).
\end{split}
\end{equation*}
Noting $\lim_{R\to \infty}\PP(\{\tau(R)\leq T\}) = 0,$ %Therefore, we replace $\tau_{n}(R)$ by $\tau_{n}(R) \wedge \tau(R)$ for $n \in %\mathbb{N}$ and $R \in [0,\infty],$ then we still have that $\mbox{\ \ for\ all\ \ } T\in[0,\infty),$
we obtain
\begin{equation}
\lim_{R\rightarrow\infty}\overline{\lim_{n\rightarrow\infty}}\PP\left(\left\{\tau_{n}(R) \wedge \tau(R)\leq T,\ \sup_{t\in[0,\tau_{n}(R) \wedge \tau(R)]}|X_{n}(t)-D(X_{n}(t-\tau))|<r(R)\right\}\right)=0.
\end{equation} 

Therefore all three assumptions hold if we replace $\tau_{n}(R)$ by $\tau_{n}(R) \wedge \tau(R).$ We may assume that $\tau_{n}(R)\leq \tau(R),$ then for $\forall t \in [0,T],$ $R \in [0,\infty)$ and $n \in \mathbb{N},$
\begin{equation}
\alpha_{R}(t\wedge\tau_{n}(R)) \leq u(R) .
\end{equation}
Now,  for fixed  $R \in [0,\infty)$, we define 
\begin{equation}\label{eq13}
\lambda_{n}^{R}(t) = \int_{0}^{t}|p_{n}(s)|K_{R}(s)ds,\quad t\in [0,\tau_{n}(R)\wedge T],\quad n\in \mathbb{N}.
\end{equation}
%For all $ T \in [0,\infty)$ and $R \in [0,\infty),$ there exists a 
Let $m \in [0,\infty),$  one has
\begin{equation}\label{eq14}
\begin{split}
&\int_{0}^{T\wedge \tau_{n}(R)}|p_{n}(s)|K_{R}(s)ds\\
&\leq m \int_{0}^{T\wedge \tau_{n}(R)}|p_{n}(s)|ds+R\int_{0}^{T\wedge \tau(R)}I_{[m,\infty)}K_{R}(s)ds.
\end{split}
\end{equation}
The right hand side of $(\ref{eq14})$ converges to 
\begin{equation}\label{eq11}
R\int_{0}^{T\wedge \tau(R)}I_{[m,\infty)}K_{R}(s)ds,
\end{equation}
 due to assumption {\bf (ii)}. Also, by observation that $(\ref{eq11})$ is dominated by $R\times u(R)$, therefore Lebesgue's dominated convergence theorem yields,
\begin{equation}\label{eq16}
\lim_{n \to \infty}\E(\lambda_{n}^{R}(T\wedge \tau_{n}(R))) = 0.
\end{equation}
%By using the same method, we denote that, for any $R\in\mathbb{R}_{+}.$ 
%\begin{equation}\label{eq13a}
%\tilde{\lambda}_{n}^{R}(t) := \int_{0}^{t}|p_{n}(s)|\tilde{K}_{R}(s)ds,\quad t\in [0,\tau_{n}(R)\wedge T],\quad n\in \mathbb{N}.
%\end{equation}
%Then, we have
%\begin{equation}\label{eq16a}
%\lim_{n \to \infty}E(\tilde{\lambda}_{n}^{R}(T\wedge \tau_{n}(R))) = 0.
%\end{equation}
%Note that, if $\tilde{K}_{R}(s)$ is replaced by $\tilde{K}_{1}(s),$ an auxiliary result can be obtained.

{\it Step (ii)}  \quad
Denote $\tau_{(n,m)}(R)=\tau_{n}(R)\wedge\tau_{m}(R),$  we now claim  that 
\begin{equation}
\sup_{t\in [0,\tau_{(n,m)}(R)\wedge T]}|X_{n}(t)-X_{m}(t)| \to 0 \quad as \quad n,m \to \infty.
\end{equation}
For simplicity, letting
\begin{equation*}
A_{m,n}(t) = X_{n}(t)-D(X_{n}(t-\tau))-X_{m}(t)+D(X_{m}(t-\tau)),
\end{equation*}
 we then have
\begin{equation*}
|X_{n}(t)-X_{m}(t)|^2= |A_{m,n}(t)+D(X_{n}(t-\tau))-D(X_{m}(t-\tau))|^2.
\end{equation*}
An application of Lemma \ref{lemma2} yields,
\begin{equation*}
|X_{n}(t)-X_{m}(t)|^2 \leq (1+\epsilon)\bigg[\frac{|D(X_{n}(t-\tau))-D(X_{m}(t-\tau))|^2}{\epsilon}+|A_{m,n}(t)|^2\bigg].
\end{equation*}
Letting $\epsilon = \frac{\kappa}{1-\kappa}$ together with assumption {\bf(C4)}, we  further obtain
\begin{equation}\label{eq18}
|X_{n}(t)-X_{m}(t)|^2 \leq \kappa|X_{n}(t-\tau)-X_{m}(t-\tau)|^2+\frac{1}{1-\kappa}|A_{m,n}(t)|^2.
\end{equation}
For a negative constant $\bar{\kappa}$, define
\begin{equation}
\psi(t) = \exp (\bar{\kappa}\alpha_{R}(t)-|\xi(0)|), \quad t \in [0,\infty).
\end{equation}
Now applying It\^{o}'s formula we have for all $t\in[0,\infty),$
\begin{equation}\label{eq17}
\begin{split}
&|A_{m,n}(t)|^2\psi(t)=\int^t_{0}\psi(s)\Big[\bar{\kappa}K_{R}(s)|A_{m,n}(s)|^2\\
&+2\langle A_{m,n}(s),b(X_{n}(s)+p_{n}(s),X_{n}(s-\tau)+p_{n}(s-\tau),s)\\
&~~~~~~~~~~~~-b(X_{m}(s)+p_{m}(s),X_{m}(s-\tau)+p_{m}(s-\tau),s)\rangle\\
&+\|\sigma(X_{n}(s)+p_{n}(s),X_{n}(s-\tau)+p_{n}(s-\tau),s)\\
&~~~~~~~~~~~-\sigma(X_{m}(s)+p_{m}(s),X_{m}(s-\tau)+p_{m}(s-\tau),s)\|^2\Big]ds\\
&+M_{n,m}^R(t),
\end{split}
\end{equation}
where 
\begin{equation*}
\begin{split}
M_{m,n}^R(t)&=\int^t_{0}2\psi(s)\langle A_{m,n}(s),\sigma(X_{n}(s)+p_{n}(s),X_{n}(s-\tau)+p_{n}(s-\tau),s)\\
&-\sigma(X_{m}(s)+p_{m}(s),X_{m}(s-\tau)+p_{m}(s-\tau),s)\rangle dB(s),
\end{split}
\end{equation*}
for $t\in[0,\infty)$ is a local $(\mathcal {F}_t)$-martingale vanishing at $t = 0,$ i.e.
$M_{n,m}^R(0)=0.$\\
Note that $A_{m,n}(s)$ can be rewritten as
\begin{equation*}
\begin{split}
&A_{m,n}(s) = X_{n}(s)-D(X_{n}(s-\tau))-X_{m}(s)+D(X_{m}(s-\tau))-p_{n}(s)+p_{m}(s)\\
&+p_{n}(s)-p_{m}(s)-D(x_{n}(s-\tau)+p_{n}(s-\tau))+D(x_{n}(s-\tau)+p_{n}(s-\tau))\\
&+D(x_{m}(s-\tau)+p_{m}(s-\tau))-D(x_{m}(s-\tau)+p_{m}(s-\tau)),
\end{split}
\end{equation*}
then by assumptions {\bf(C3)}, {\bf(C4)}, we have for any $t\in[0,\tau_{(n,m)}(R)\wedge T],$
\begin{equation}\label{eq15}
\begin{split}
&|A_{m,n}(t)|^2\psi(t) = \int^t_{0}\psi(s)\Big[\bar{\kappa}K_{R}(s)|A_{m,n}(s)|^2\\
&+2\langle-p_{n}(s)+p_{m}(s),b(X_{n}(s)+p_{n}(s),X_{n}(s-\tau)+p_{n}(s-\tau),s)\\
&-b(X_{m}(s)+p_{m}(s),X_{m}(s-\tau)+p_{m}(s-\tau),s)\rangle\\
&+2\langle D(X_{m}(s-\tau))-D(X_{m}(s-\tau)+p_{m}(s-\tau)),b(X_{n}(s)+p_{n}(s),X_{n}(s-\tau)+p_{n}(s-\tau),s)\\
&-b(X_{m}(s)+p_{m}(s), X_{m}(s-\tau)+p_{m}(s-\tau),s)\rangle\\
&+2\langle D(X_{n}(s-\tau)+p_{n}(s-\tau))-D(X_{n}(s-\tau)),\\
&  b(X_{n}(s)+p_{n}(s),X_{n}(s-\tau)+p_{n}(s-\tau),s)-b(X_{m}(s)+p_{m}(s),X_{m}(s-\tau)+p_{m}(s-\tau),s)\rangle\\
&+2\langle X_{n}(s)+p_{n}(s)-D(X_{n}(s-\tau)+p_{n}(s-\tau))-X_{m}(s)-p_{m}(s)+D(X_{m}(s-\tau)+p_{m}(s-\tau)),\\
& b(X_{n}(s)+p_{n}(s),X_{n}(s-\tau)+p_{n}(s-\tau),s)-b(X_{m}(s)+p_{m}(s),X_{m}(s-\tau)+p_{m}(s-\tau),s) \rangle\\
&+\|\sigma(X_{n}(s)+p_{n}(s),X_{n}(s-\tau)+p_{n}(s-\tau),s)\\
&-\sigma(X_{m}(s)+p_{m}(s),X_{m}(s-\tau)+p_{m}(s-\tau),s)\|^{2}\Big]ds+M_{n,m}^R(t)\\
&\leq \int^t_{0}\psi(s)\Big[\bar{\kappa}K_{R}(s)|A_{m,n}(s)|^2\\
&+2\langle-p_{n}(s)+p_{m}(s),b(X_{n}(s)+p_{n}(s),X_{n}(s-\tau)+p_{n}(s-\tau),s)\\
&-b(X_{m}(s)+p_{m}(s),X_{m}(s-\tau)+p_{m}(s-\tau),s)\rangle\\
&+2\langle D(X_{m}(s-\tau))-D(X_{m}(s-\tau)+p_{m}(s-\tau)),b(X_{n}(s)+p_{n}(s),X_{n}(s-\tau)+p_{n}(s-\tau),s)\\
&-b(X_{m}(s)+p_{m}(s), X_{m}(s-\tau)+p_{m}(s-\tau),s)\rangle\\
&+2\langle D(X_{n}(s-\tau)+p_{n}(s-\tau))-D(X_{n}(s-\tau)),\\
&  b(X_{n}(s)+p_{n}(s),X_{n}(s-\tau)+p_{n}(s-\tau),s)-b(X_{m}(s)+p_{m}(s),X_{m}(s-\tau)+p_{m}(s-\tau),s)\rangle\\
&+K_{R}(s)|X_{n}(s)+p_{n}(s)-X_{m}(s)-p_{m}(s)|^{2}\\
&+\tilde{K}_{R}(s-\tau)|X_{n}(s-\tau)+p_{n}(s-\tau)-X_{m}(s-\tau)-p_{m}(s-\tau)|^{2}\Big]ds+M_{n,m}^R(t)\\
%&\leq \int^t_{0}\psi(s)\Big[\bar{\kappa}K_{R}(s)|A_{m,n}(s)|^2+4K_{R}(s)|p_{m}(s)-p_{n}(s)|\\
%&+4 K_{R}(s)|D(X_{m}(s-\tau))-D(X_{n}(s-\tau))-(D(X_{m}(s-\tau)+p_{m}(s-\tau)\\
%&-D(X_{n}(s-\tau)+p_{n}(s-\tau)))|+K^{(1)}_{R}(s)|X_{n}(s)+p_{n}(s)-X_{m}(s)-p_{m}(s)|^{2}\\
%&+K_{R}(s)|X_{n}(s-\tau)+p_{n}(s-\tau)-X_{m}(s-\tau)+p_{m}(s-\tau)|^{2}\Big]ds+M_{n,m}^R(t)\\
&\leq\int^t_{0}\psi(s)\Big[\bar{\kappa}K_{R}(s)|A_{m,n}(s)|^2+4K_{R}(s)(|p_{m}(s)-p_{n}(s)|\\
&+\kappa|p_{m}(s-\tau)|+\kappa|p_{n}(s-\tau)|)+2K_{R}(s)(|X_{n}(s)-X_{m}(s)|^2+|p_{n}(s)-p_{m}(s)|^{2})\\
&+2\tilde{K}_{R}(s-\tau)(|X_{n}(s-\tau)-X_{m}(s-\tau)|^2+|p_{n}(s-\tau)-p_{m}(s-\tau)|^2)\Big]ds+M_{n,m}^R(t).\end{split}
\end{equation}
By \eqref{eq18}, we derive that 
\begin{equation}\label{eq31}
\begin{split}
& |A_{m,n}(t)|^2\psi(t)\\
&\leq  \int^t_{0}\psi(s)\Big[\bar{\kappa}K_{R}(s)\bigg((1-\kappa)|X_{n}(s)-X_{m}(s)|^2-\kappa(1-\kappa) |X_{n}(s-\tau)-X_{m}(s-\tau)|^2 \bigg)\\
&+4K_{R}(s)(|p_{m}(s)-p_{n}(s)|+\kappa|p_{m}(s-\tau)|+\kappa|p_{n}(s-\tau)|)\\
&+2K_{R}(s)(|X_{n}(s)-X_{m}(s)|^2+|p_{n}(s)-p_{m}(s)|^{2})\\
&+2\tilde{K}_{R}(s-\tau)(|X_{n}(s-\tau)-X_{m}(s-\tau)|^2+|p_{n}(s-\tau)-p_{m}(s-\tau)|^2)\Big]ds+M_{n,m}^R(t),\end{split}
\end{equation}
%For the following term:
%\begin{equation*}
%\int_{0}^{t}\psi(s)K_{R}(s)\bar{\kappa}(\kappa^{2}-\kappa)|X_{n}(s-\tau)-X_{m}(s-\tau)|^2ds
%\end{equation*}
Since for any $s\in[0,t],$ $\psi(s)$ is a non-increasing function, also note that $X_n(t)\equiv X_m(t),$ for any $ t \in [-\tau, 0],$  we have 
\begin{equation}\label{eq21}
\begin{split}
&\int_{0}^{t}\psi(s)K_{R}(s)\bar{\kappa}(\kappa^{2}-\kappa)|X_{n}(s-\tau)-X_{m}(s-\tau)|^2ds\\
&\leq C_{R}(\tau)\int_{0}^{t}\psi(s-\tau)K_{R}(s-\tau)\bar{\kappa}(\kappa^{2}-\kappa)|X_{n}(s-\tau)-X_{m}(s-\tau)|^2ds\\
&\leq C_{R}(\tau) \int_{0}^{t}\psi(s)K_{R}(s)\bar{\kappa}(\kappa^{2}-\kappa)|X_{n}(s)-X_{m}(s)|^2ds,
\end{split}
\end{equation} 
%Now, apply the same technique to the term
%\begin{equation*}
%\int_{0}^{t}\psi(s)K_{R}(s)|X_{n}(s-\tau)-X_{m}(s-\tau)|^2ds
%\end{equation*}
%we have
and 
\begin{equation}\label{eq30}
\begin{split}
&\int_{0}^{t}\psi(s)\tilde{K}_{R}(s-\tau)(|X_{n}(s-\tau)-X_{m}(s-\tau)|^2+|p_{n}(s-\tau)-p_{m}(s-\tau)|^2)ds\\
&\le\int_{0}^{t}\psi(s-\tau)\tilde{K}_{R}(s-\tau)(|X_{n}(s-\tau)-X_{m}(s-\tau)|^2+|p_{n}(s-\tau)-p_{m}(s-\tau)|^2)ds\\
&\leq 
\int_{0}^{t}\psi(s)\tilde{K}_{R}(s)(|X_{n}(s)-X_{m}(s)|^2+|p_{n}(s)-p_{m}(s)|^{2})ds.
\end{split}
\end{equation}
%since $K_{R}(s)\geq K^{(1)}_{R}(s)$ for all $s\in [0,T],$ therefore, 
%\begin{equation*}
%2K^{(1)}_{R}(s)(|X_{n}(s)-X_{m}(s)|^2+|p_{n}(s)-p_{m}(s)|^{2})\leq 2K_{R}(s)(|X_{n}(s)-X_{m}(s)|^2+|p_{n}(s)-p_{m}(s)|^{2})
%\end{equation*}

Now substituting $(\ref{eq21})$ and $(\ref{eq30})$  into $(\ref{eq31}),$ which yields
\begin{equation*}
\begin{split}
|A_{m,n}(t)|^2\psi(t)&\leq  \int^t_{0}\psi(s)\bigg[K_{R}(s)(\bar{\kappa}((1-\kappa)+C_{R}(\tau)(\kappa^2-\kappa))+2)|X_{n}(s)-X_{m}(s)|^2\\
&+4K_{R}(s)(|p_{m}(s)-p_{n}(s)|+\kappa |p_{m}(s-\tau)|+\kappa |p_{n}(s-\tau)|)\\
&+2K_{R}(s)|p_{n}(s)-p_{m}(s)|^{2}\\
&+2\tilde{K}_{R}(s)(|X_{n}(s)-X_{m}(s)|^2+|p_{n}(s)-p_{m}(s)|^{2})\bigg]ds
+M_{n,m}^R(t).
\end{split}
\end{equation*}
Noting that $K_{R}(t)\geq \tilde{K}_{R}(t)$ and  choosing $\bar{\kappa} = \frac{-4}{((1-\kappa)+C_{R}(\tau)(\kappa^2-\kappa))},$ we obtain
\begin{equation*}
\begin{split}
|A_{m,n}(t)|^2\psi(t)&\leq  \int^t_{0}4K_{R}(s)\psi(s)(|p_{m}(s)-p_{n}(s)|+\kappa |p_{m}(s-\tau)|+\kappa |p_{n}(s-\tau)|\\
&+|p_{n}(s)-p_{m}(s)|^{2})ds+M_{n,m}^R(t).
\end{split}
\end{equation*}
It is easy to see
\begin{equation*}
\begin{split}
&\int^t_{0}4\kappa\psi(s)K_{R}(s)(|p_{n}(s-\tau)|+|p_{m}(s-\tau)|)\\
&\leq \int^t_{0}C_{R}(\tau)\psi(s-\tau)4kK_{R}(s-\tau)(|p_{n}(s-\tau)|+|p_{m}(s-\tau)|)\\
&\leq 4\kappa C_{R}(\tau)\int^t_{0}\psi(s)K_{R}(s)(|p_{n}(s)|+|p_{m}(s)|).
\end{split}
\end{equation*}
and for all $ t \in [0,\tau_{(n,m)}(R)\wedge T],$ $\psi(t)< 1$ 
\begin{equation*}
|-p_{n}(t)+p_{m}(t)|^{2}\leq 2R(|p_{n}(t)|+|p_{m}(t)|) \quad a.s..
\end{equation*}
Then we  have for $t \in [0, \tau_{(n,m)}(R)\wedge T],$
\begin{equation*}
\begin{split}
&|A_{m,n}(t)|^2\psi(t)\leq (4\kappa C_{R}(\tau)+8R+4)(\lambda_{n}(t)+\lambda_{m}(t)).
\end{split}
\end{equation*} 
  Hence for any $\mathcal{F}_{t}-$stopping time $\bar{\tau} \leq \tau_{(n,m)}(R)$ and $\mathcal{F}_{t}-$stopping times $\sigma_{l}\uparrow \infty $ as $l \to \infty$,  $M_{n,m}^R(t\wedge \sigma_{l})$ is martingale for all $l \in \mathbb{N}.$ Therefore, we have
\begin{equation*}
\begin{split}
\E(|A_{m,n}(\bar{\tau} \wedge \sigma_{l})|^2\psi(\bar{\tau}\wedge \sigma_{l}))&\leq(4\kappa C_{R}(\tau)+8R+4) \E(\lambda_{n}(T \wedge \tau_{n}(R) )+\lambda_{m}(T\wedge \tau_{m}(R) )).\\
\end{split}
\end{equation*}
%First letting $l \to \infty,$
%\begin{equation*}
%|A_{m,n}(t)|^2\psi(t) = \inf{\lim_{l\to \infty}}|A_{m,n}(\tau \wedge \sigma_{l})|^2\psi(\tau \wedge \sigma_{l})
%\end{equation*}
Then the Fatou Lemma yields
\begin{equation*}
\begin{split}
\E|A_{m,n}(\bar{\tau})|^2\psi(\bar{\tau})&\leq{\lim_{l\to \infty}}\inf\big( \E(|A_{m,n}(\bar{\tau} \wedge \sigma_{l})|^2\psi(\bar{\tau} \wedge \sigma_{l}))\\
&\leq (4\kappa C_{R}(\tau)+8R+4)\E(\lambda_{n}(T \wedge \tau_{n}(R))+\lambda_{m}(T\wedge \tau_{m}(R))).
\end{split}
\end{equation*}
Then using Lemma \ref{lemma1}, we obtain that for every $\epsilon\in (0,\infty)$
\begin{equation}\label{eq19}
\begin{split}
& \PP({\sup_{t\in [0,\tau_{(n,m)}(R)\wedge T]}(|A_{m,n}(t)|^2\psi(t)) > \epsilon})\\
& = \frac{1}{\epsilon}[(4\kappa C_{R}(\tau)+8R+4)\E(\lambda_{n}(T \wedge \tau_{n}(R))+\lambda_{m}(T \wedge \tau_{m}(R))).
\end{split}
\end{equation}
Since $[0,T] \ni t \mapsto \psi(t)$ is strictly positive, which is independent of $n,m \in \mathbb{N}$ and continuous,  $(\ref{eq19})$ implies that 
\begin{equation*}
\sup_{t\in [0,\tau_{(n,m)}(R)\wedge T]}|A_{m,n}(t)| \xrightarrow{\mathbb{P}} 0 \quad n,m \to \infty.
\end{equation*}
Recall that $X_{n}(t) \equiv X_{m}(t)$ for $t\in [-\tau,0],$ based on the fact given by Lemma \ref{LemmaKey}, we have 
\begin{equation*}
\sup_{t\in [0,\tau_{(n,m)}(R)\wedge T]}|X_{n}(t)-X_{m}(t)| \xrightarrow{\mathbb{P}} 0 \quad n,m \to \infty.
\end{equation*}

{\it Step (iii)}\quad 
 we shall show that for any given $T \in [0,\infty),$
\begin{equation*}
\lim_{R \to \infty}\overline{\lim}_{n \to \infty}\PP({\tau_{n}(R)\leq T}) = 0.
\end{equation*}
By using Lemma $\ref{lemma2},$ we have 
\begin{equation}
\begin{split}
|X_n(t)|^2\leq (1+\epsilon)\left(\frac{|D(X_{n}(t-\tau))|^2}{\epsilon}+|X_{n}(t)-D(X_{n}(t-\tau))|^2\right).
\end{split}
\end{equation}
Letting $\epsilon = \frac{\kappa}{1-\kappa} $ and noting the assumption {\bf(C4)}, 
we derive
\begin{equation}\label{eq24}
|X_{n}(t)-D(X_{n}(t-\tau))|^2\geq  (\kappa^{2}-\kappa) |X_{n}(t-\tau)|^{2}+(1-\kappa)|X_{n}(t)|^{2}.
\end{equation}
%that is
%\begin{equation}
%\sup_{t\in[0,T]}|X_{n}(t)|^{2}\leq  \frac{\kappa}{1-\kappa}\sup_{t\in[-\tau,0]}|\xi(t)|^{2}+\bigg(\frac{1}{1-\kappa}\bigg)^2\sup_{t\in[0,T]}|X_{n}(t)-D(X_{n}(t-\tau))|^2.
%\end{equation}
 Let $\tilde{\kappa}$ be a negative constant and define \begin{equation*}
\varphi(t) = \exp (\tilde{\kappa}\alpha_{1}(t)-|\xi(0)|), \quad t \in [0,\infty),
\end{equation*}
where $\alpha_{1}(t) = \int_{0}^{t}K_{1}(s)ds.$
 An application of  It\^{o}'s formula implies
 \begin{equation}\label{eq27}
\begin{split}
&|X_{n}(t)-D(X_{n}(t-\tau))|^2\varphi(t) = |\xi(0)-D(\xi(-\tau))|^2e^{-|\xi(0)|}\\
&+\int_{0}^{t}\varphi(s)[2\langle X_{n}(s)-D(X_{n}(s-\tau)),b(X_{n}(s)+p_{n}(s),X_{n}(s-\tau)+p_{n}(s-\tau),s)\rangle\\
&+||\sigma(X_{n}(s)+p_{n}(s),X_{n}(s-\tau)+p_{n}(s-\tau),s)||^2\\
&+\tilde{\kappa}K_{1}(t)|X_{n}(s)-D(X_{n}(s-\tau))|^2]ds+M_{n}^{R}(t),
\end{split}
\end{equation}
where %$M_{n}^{R}(t):$
\begin{equation*}
\begin{split}
M_{n}^R(t)&=\int^t_{0}2\varphi(s)\langle X_{n}(s)-D(X_{n}(s-\tau)),\sigma(X_{n}(s)+p_{n}(s),X_{n}(s-\tau)+p_{n}(s-\tau),s)\rangle dB(s),
\end{split}
\end{equation*}
for $t\in[0,\infty)$ is a local $\mathcal{F}_{t}$-martingale with $M_{n}^{R}(0)= 0$. Using assumption {\bf(C2)} and hypothesis {\bf(H)}  for all $t \in [0,T\wedge \tau_{n}(R)]$, we compute 
\begin{equation}\label{eq20}
\begin{split}
&|X_{n}(t)-D(X_{n}(t-\tau))|^2\varphi(t)\\
&=|\xi(0)-D(\xi(-\tau))|^2e^{-|\xi(0)|}+\int^t_{0}\varphi(s)\Big[\tilde{\kappa}K_{1}(s)|X_{n}(s)-D(X_{n}(s-\tau))|^2\\
&+2\langle X_{n}(s)-D(X_{n}(s-\tau))+p_{n}(s)-p_{n}(s)-D(X_{n}(s-\tau)\\
&+p_{n}(s-\tau))+D(X_{n}(s-\tau)+p_{n}(s-\tau)),
b(X_{n}(s)+p_{n}(s),X_{n}(s-\tau)+p_{n}(s-\tau),s)\rangle\\
&+\|\sigma(X_{n}(s)+p_{n}(s),X_{n}(s-\tau)+p_{n}(s-\tau),s)\|^2\Big]ds+M_{n}^R(t)\\
&= |\xi(0)-D(\xi(-\tau))|^2e^{-|\xi(0)|}+\int^t_{0}\varphi(s)\Big[\tilde{\kappa}K_{1}(s)|X_{n}(s)-D(X_{n}(s-\tau))|^2\\
&+2\langle X_{n}(s)+p_{n}(s)-D(X_{n}(s-\tau)+p_{n}(s-\tau)),b(X_{n}(s)+p_{n}(s),X_{n}(s-\tau)+p_{n}(s-\tau),s)\rangle\\
&+2\langle -D(X_{n}(s-\tau))-p_{n}(s)+D(X_{n}(s-\tau)+p_{n}(s-\tau)),\\
& b(X_{n}(s)+p_{n}(s),X_{n}(s-\tau)+p_{n}(s-\tau),s)\rangle \\
&+\|\sigma(X_{n}(s)+p_{n}(s),X_{n}(s-\tau)+p_{n}(s-\tau),s)\|^2\Big]ds+M_{n}^R(t)\\
& \leq |\xi(0)-D(\xi(-\tau))|^2e^{-|\xi(0)|}+\int^t_{0}\varphi(s)\Big[\tilde{\kappa}K_{1}(s)|X_{n}(s)-D(X_{n}(s-\tau))|^2\\
& 2\langle -p_{n}(s)+D(X_{n}(s-\tau)+p_{n}(s-\tau))-D(X_{n}(s-\tau)),\\
& b(X_{n}(s)+p_{n}(s),X_{n}(s-\tau)+p_{n}(s-\tau),s)\rangle\\
&+K_{1}(s)(1+|X_{n}(s)+p_{n}(s)|^2)+\tilde{K}_{1}(s-\tau)(1+|X_{n}(s-\tau)+p_{n}(s-\tau)|^{2})\Big]ds+M_{n}^R(t).\\
%&\leq |\xi(0)-D(\xi(-\tau))|^2e^{-|\xi(0)|}+\int^t_{0}\varphi(s)\Big[\tilde{\kappa}K_{1}(s)|X_{n}(s)-D(X_{n}(s-\tau))|^2\\
%&+2K_{1}(s)|-p_{n}(s)+D(X_{n}(s-\tau)+p_{n}(s-\tau))-D(X_{n}(s-\tau))|\\
%&+K_{1}(s)(1+|X_{n}(s)+p_{n}(s)|^2)+\tilde{K}_{1}(s-\tau)(1+|X_{n}(s-\tau)+p_{n}(s-\tau)|^{2})\Big]ds+M_{n}^R(t).
\end{split}
\end{equation}
%This, together with $(\ref{eq24}),$ implies
%\begin{equation}\label{eq28}
%\begin{split}
%&|X_{n}(t)-D(X_{n}(t-\tau))|^2\varphi(t)\\
%&\le|X_{0}-D(\xi)|^2e^{-|X_{0}|}+\int^t_{0}\varphi(s)\Big[\tilde{\kappa}K_{1}(s)((\kappa^{2}-\kappa) |X_{n}(s-\tau)|^{2}+(1-\kappa)|X_{n}(s)|^{2})\\
%&+2K^{(1)}_{1}(s)|-p_{n}(s)+D(X_{n}(s-\tau)+p_{n}(s-\tau))-D(X_{n}(s-\tau))|\\
%&+K^{(1)}_{1}(s)(1+|X_{n}(s)+p_{n}(s)|^2)+K_{1}(s)(1+|X_{n}(s-\tau)+p_{n}(s-\tau)|^{2})\Big]ds+M_{n}^R(t)\\
%\end{split}
%\end{equation}
 Using  \eqref{eq24}, we can write that for $t\in[0,\tau_{n}(t)\wedge T]$
%together with the fact that $K_{1}(s)\geq K^{(1)}_{1}(s)$ for all $s\in [0,T],$ therefore, 
%\begin{equation*}
%K^{(1)}_{1}(s)(1+|X_{n}(s)+p_{n}(s)|^2)\leq K_{1}(s)(1+|X_{n}(s)+p_{n}(s)|^2)
%\end{equation*} 
\begin{equation}\label{eq29}
\begin{split}
&|X_{n}(t)-D(X_{n}(t-\tau))|^2\varphi(t)\\
&\leq |\xi(0)-D(\xi(-\tau))|^2e^{-|\xi(0)|}+\int^t_{0}\varphi(s)\Big[\tilde{\kappa}K_{1}(s)((\kappa^{2}-\kappa) |X_{n}(s-\tau)|^{2}+(1-\kappa)|X_{n}(s)|^{2})\\
&+2K_{R}(s)|-p_{n}(s)+D(X_{n}(s-\tau)+p_{n}(s-\tau))-D(X_{n}(s-\tau))|\\
&+K_{1}(s)(1+|X_{n}(s)+p_{n}(s)|^2)+\tilde{K}_{1}(s-\tau)(1+|X_{n}(s-\tau)+p_{n}(s-\tau)|^{2})\Big]ds+M_{n}^R(t).
\end{split}
\end{equation}
Recalling that $\varphi(t)$ is non-increasing for all $t\in[0,\infty),$  we can write that
\begin{equation*}
\begin{split}
&\int_{0}^{t}\varphi(s)\tilde{\kappa}K_{1}(s)(\kappa^{2}-\kappa) |X_{n}(s-\tau)|^{2}ds\\
&\leq C_{1}(\tau) \int_{0}^{t}\varphi(s-\tau)\tilde{\kappa}K_{1}(s-\tau)(\kappa^{2}-\kappa) |X_{n}(s-\tau)|^{2}ds\\
%& = C_{1}(\tau) \int_{0}^{t-\tau}\varphi(s)\tilde{\kappa}K_{1}(s)(\kappa^{2}-\kappa) |X_{n}(s)|^{2}ds\\
%&+ C_{1}(\tau) \int_{-\tau}^{0}\varphi(\theta)\tilde{\kappa}K_{1}(\theta)(\kappa^{2}-\kappa) |\xi(\theta)|^{2}d\theta\\
&\leq  C_{1}(\tau) \int_{0}^{t}\varphi(s)\tilde{\kappa}K_{1}(s)(\kappa^{2}-\kappa) |X_{n}(s)|^{2}ds
+ C_{1}(\tau) \int_{-\tau}^{0}\varphi(\theta)\tilde{\kappa}K_{1}(\theta)(\kappa^{2}-\kappa) |\xi(\theta)|^{2}d\theta,
\end{split}
\end{equation*}
\begin{equation*}
\begin{split}
&\int_{0}^{t}\varphi(s)\tilde{K}_{1}(s-\tau)(1+|X_{n}(s-\tau)+p_{n}(s-\tau)|^{2})ds\\
%&\leq \int_{0}^{t}\varphi(s-\tau)\tilde{K}_{1}(s-\tau)(1+|X_{n}(s-\tau)+p_{n}(s-\tau)|^{2})ds\\
%&\leq \int_{0}^{t-\tau}\varphi(s)\tilde{K}_{1}(s)(1+|X_{n}(s)+p_{n}(s)|^{2})ds\\
%&+\int_{-\tau}^{0}\varphi(\theta)\tilde{K}_{1}(\theta)(1+|\xi(\theta)|^{2})d\theta\\
&\leq \int_{0}^{t}\varphi(s)\tilde{K}_{1}(s)(1+|X_{n}(s)+p_{n}(s)|^{2})ds
+\int_{-\tau}^{0}\varphi(\theta)\tilde{K}_{1}(\theta)(1+|\xi(\theta)|^{2})d\theta.
\end{split}
\end{equation*}
\begin{equation*}
\int_{0}^{t}2\kappa\varphi(s)K_{R}(s)|p_{n}(s-\tau)|ds\leq C_{R}(\tau)\int_{0}^{t}2\kappa\varphi(s)K_{R}(s)|p_{n}(s)|ds.
\end{equation*}
Therefore, we can rewrite \eqref{eq29} as for $t\in[0,\tau_{n}(t)\wedge T]$
\begin{equation}\label{eq32}
\begin{split}
&|X_{n}(t)-D(X_{n}(t-\tau))|^2\varphi(t)\\
&\leq |\xi(0)-D(\xi(-\tau))|^2e^{-|\xi(0)|}+\int_{0}^{t}\varphi(s)\bigg[K_{1}(s)(\tilde{\kappa}((1-\kappa)+C_{1}(\tau)(\kappa^2-\kappa))+2)|X_{n}(t)|^{2}\\ &+2K_{R}(s)(|p_{n}(s)|+\kappa C_R(\tau)|p_{n}(s)|)+K_{1}(s)(1+2|p_{n}(s)|^{2})\bigg]ds\\
&+ C_{1}(\tau)\int_{-\tau}^{0}\varphi(\theta)\tilde{\kappa}K_{1}(\theta)(\kappa^{2}-\kappa) |\xi(\theta)|^{2}d\theta +\int_{0}^{t}\varphi(s)\tilde{K}_{1}(s)(1+2|X_{n}(s)|^{2}+2|p_{n}(s)|^{2})ds\\
&+\int_{-\tau}^{0}\varphi(\theta)\tilde{K}_{1}(\theta)(1+|\xi(\theta)|^{2})d\theta+M_{n}^R(t).
\end{split}\end{equation}
%Recall that
%\begin{equation*}
%\int_{0}^{t}2\kappa\varphi(s)K_{R}(s)|p_{n}(s-\tau)|ds\leq C_{R}(\tau)\int_{0}^{t}2\kappa\varphi(s)K_{R}(s)|p_{n}(s)|ds.
%\end{equation*}

Again, since $T$ is fixed, then for any $t\in[0,T],$ noting that $K_{1}(t)\geq \tilde{K}_{1}(t),$ then by choosing $\bar{\kappa} =-\frac{4}{(1-\kappa)+C_{1}(\tau)(\kappa^2-\kappa)},$ we have  for $t\in[0,\tau_{n}(t)\wedge T]$
\begin{equation}\label{eq34}
\begin{split}
&|X_{n}(t)-D(X_{n}(t-\tau))|^2\varphi(t)\leq |\xi(0)-D(\xi(-\tau))|^2e^{-|\xi(0)|}+\int_{0}^{t}\varphi(s)\bigg[2K_{R}(s)(1+\kappa C_R(\tau))|p_{n}(s)|\\
&+K_{1}(s)(1+2|p_{n}(s)|^{2})\bigg]ds+ C_{1}(\tau)\int_{-\tau}^{0}\varphi(\theta)\tilde{\kappa}K_{1}(\theta)(\kappa^{2}-\kappa) |\xi(\theta)|^{2}d\theta\\ &+\int_{0}^{t}\varphi(s)\tilde{K}_{1}(s)(1+2|p_{n}(s)|^{2})ds+\int_{-\tau}^{0}\varphi(\theta)\tilde{K}_{1}(\theta)(1+|\xi(\theta)|^{2})d\theta+M_{n}^R(t).
\end{split}
\end{equation}
Then for $t \in [0,T],$ without losing generality, we may replace 
$K_{R}(t)$ by the $\max\{K_{R}(t),K_{1}(t),\tilde{K}_{1}(t)\},$ then we can deduce that for every $\mathcal{F}_{t}-$stopping time $\tilde{\tau}\leq T\wedge\tau_{n}(R),$
\begin{equation*}
\begin{split}
& \E|X_{n}(\tilde{\tau})-D(X_{n}(\tilde{\tau}-\tau))|^2\varphi(\tau)\leq \E|\xi(0)-D(\xi(-\tau))|^2e^{-|\xi(0)|}\\
&+(2C_{R}(\tau)\kappa+2+4R)\E(\lambda_{n}^{R}(T\wedge\tau_{n}(R)))+\int_{0}^{t}2\varphi(s)K_{R}(s)ds.\\
&+\E\int_{-\tau}^{0}\varphi(\theta)K_{R}(\theta)(1+|\xi(\theta)|^{2})d\theta+C_{1}(\tau)\E \int_{-\tau}^{0}\varphi(\theta)K_{R}(\theta)(\kappa^{2}-\kappa)\tilde{\kappa}|\xi(\theta)|^{2}d\theta.\\
\end{split}
\end{equation*}
%where $C = C(\tilde{\kappa},\kappa,E|\xi(\theta)|^{2})$ is a general constant.
Therefore, by using Lemma \ref{lemma1} and \eqref{eq16}, we obtain that $\forall c \in (0,\infty),$
\begin{equation*}
\lim_{c \to \infty}\sup_{R\in [0,\infty)}\overline{\lim}_{n \to \infty}\PP({\sup_{t\in [0,T\wedge\tau_{n}(R)]}(X_{n}(t)-D(X_{n}(t-\tau))|^2\varphi(t))\geq c}) = 0.
\end{equation*}
Since $[0,T\wedge\tau_{n}(R)]\ni t \mapsto \varphi(t)$ is strictly positive and it is independent of $n \in \mathbb{N}$ and continuous, also we recall that $r(R) \to \infty$ as $R \to \infty,$ we conclude that 
\begin{equation*}
\begin{split}
&\overline{\lim}_{R\to \infty}\overline{\lim}_{n \to \infty}\PP({\sup_{t\in [0,\tau_{n}(R)]}|X_{n}(t)-D(X_{n}(t-\tau))|\geq r(R),\tau_{n}(R) \leq T})\\
& \leq \lim_{R \to \infty} \sup_{\bar{R}\in [0,\infty]}\overline{\lim}_{k \to \infty}\PP({\sup_{t\in [0,T\wedge\tau_{n}(\bar{R})]}(|X_{n}(t)-D(X_{n}(t-\tau))|^2)\geq r(R)}) = 0.
\end{split}
\end{equation*}
Therefore, by assumption {\bf (iii)} we have shown that
\begin{equation}\label{eq22}
\lim_{R \to \infty}\overline{\lim}_{n \to \infty}\PP({\tau_{n}(R)\leq T}) = 0.
\end{equation}
Hence, we complete the proof of the localization lemma. $\Box$

\section{Proof of Existence and Uniqueness theorem}

Having the localization lemma in hand, we can now prove Theorem \ref{EUtheorem}.

{\bf Proof of Theorem \ref{EUtheorem}}. Let $T>0$ be fixed such that $T/\tau$ is a rational number. Let the step size $\Delta\in (0, 1)$ be fraction of $\tau$ and $T,$ i.e. there exist positive integer $M, N$ such that $\Delta=T/M=\tau/N.$
The discrete-time Euler scheme is defined as follows:
\begin{equation}\label{yeu}
\begin{cases}
&X^{\Delta}(t) = \xi(t),-\tau\leq t\leq 0,\\
&X^{\Delta}(((l+1)\Delta)) =D(X^{\Delta}((l+1)\Delta-\tau)) + X^{\Delta}(l\Delta)-D(X^{\Delta}(l\Delta-\tau))\\
&+b(X^{\Delta}(l\Delta),X^{\Delta}(l\Delta-\tau),l\Delta)\Delta
+\sigma(X^{\Delta}(l\Delta),X^{\Delta}(l\Delta-\tau),l\Delta)\Delta B_l, 0\le l \le M,
\end{cases}
\end{equation}
where $\D B_l=B((l+1)\D)-B(l\D).$   For $a>0,$ let $[a]$ be the integer part of $a.$ Define  $\kappa(\Delta,t)= \big[\frac{t}{\Delta}\big]\Delta.$ Then we can define the continuous-time approximation of equation \eqref{yeu} as follows:
\begin{equation}\label{eq23}
\begin{split}
X^{\Delta}(t) &=D(X^{\Delta}(t-\tau)) + \xi(0)-D(\xi(-\tau))\\
&+ \int_{0}^{t}b(X^{\Delta}(\kappa(\Delta,s)),X^{\Delta}(\kappa(\Delta,s)-\tau),s)ds\\
&+\int_{0}^{t}\sigma(X^{\Delta}(\kappa(\Delta,s)),X^{\Delta}(\kappa(\Delta,s)-\tau),s)dB(s).
\end{split}
\end{equation}
Then for any $t\in[-\tau,0),$ define $p^{\Delta}(t) = 0$  and for $t \in [0,T],$  define
\begin{equation*}
\begin{split}
p^{\Delta}(t) = X^{\Delta}(\kappa(\Delta,t))-X^{\Delta}(t). 
\end{split}
\end{equation*}
As a result, \eqref{eq23} is equivalent to 
\begin{equation}\label{eq233}
\begin{split}
X^{\Delta}(t) &=D(X^{\Delta}(t-\tau)) + \xi(0)-D(\xi(-\tau))\\
&+ \int_{0}^{t}b(X^{\Delta}(s)+p^{\Delta}(s),X^{\Delta}(s-\tau)+p^{\Delta}(s-\tau,s))ds\\
&+\int_{0}^{t}\sigma(X^{\Delta}(s)+p^{\Delta}(s),X^{\Delta}(s-\tau)+p^{\Delta}(s-\tau),s)dB(s).
\end{split}
\end{equation}
Note that
\begin{equation*}
\begin{split}
p^{\Delta}(t) &= X^{\Delta}(\kappa(\Delta,t))-X^{\Delta}(t)\\
&=D(X^{\Delta}(\kappa(\Delta,t))-\tau)-D(X^{\Delta}(t)-\tau)\\ 
&-\int_{\kappa(\Delta,t)}^{t}b(X^{\Delta}(\kappa(\Delta,s)),X^{\Delta}(\kappa(\Delta,s)-\tau),s)ds\\
&-\int_{\kappa(\Delta,t)}^{t}\sigma(X^{\Delta}(\kappa(\Delta,s)),X^{\Delta}(\kappa(\Delta,s)-\tau),s)dB(s).
\end{split}
\end{equation*}
Fix $R \in [0,\infty)$, and define that
\begin{equation*}
\tau^{\Delta}(R) := \inf\{t \geq 0 | |X^{\Delta}(t)| > \frac{R}{3}\}.
\end{equation*}
Then clearly, for $\forall t \in [0,\tau^{\Delta}(R)]$
\begin{equation*}
|X^{\Delta}(t)| \leq \frac{R}{3} \quad \mbox{and} \quad |p^{\Delta}(t)| \leq \frac{2R}{3}.
\end{equation*}
As a result of that, assumption {\bf(i)} in the localization lemma holds. We may assume that $||\xi||_{(-\tau,0)}\le R/3,$ and set $r(R)$ as the following function,
\begin{equation*}
r(R):= (1-\kappa)\frac{R}{3}.
\end{equation*}
Since 
$$
\sup_{t\in [0, \tau_n(R)]}|X^{\Delta}(t)-D(X^{\Delta}(t-\tau))|\ge \sup_{t\in [0, \tau_n(R)]}|X^{\Delta}(t)|-\kappa\sup_{t\in [0, \tau_n(R)]} |X^{\Delta}(t-\tau)|\ge r(R),
$$
 the assumption {\bf(iii)} in the localization lemma is empty for all $\Delta \in (0,1)$ and $R \in[0,\infty),$  this means  the assumption  {\bf(iii)} is also fulfilled.

In order to   show the assumption {\bf(ii)} in the localization lemma holds, we compute
\begin{equation}\label{SecondConditionCheck}
\begin{split}
& \E\int_{0}^{T\wedge\tau^{\Delta}(R)}|p^{\Delta}(s)|ds = \E\int_{0}^{T\wedge\tau^{\Delta}(R)}|X^{\Delta}(\kappa(\Delta,s))-X^{\Delta}(s)|ds\\
&\leq \E\int_{0}^{T\wedge\tau^{\Delta}(R)}\bigg|D(X^{\Delta}(\kappa(\Delta,s))-\tau)-D(X^{\Delta}(s)-\tau)\bigg|ds\\  &+\E\int_{0}^{T\wedge\tau^{\Delta}(R)}\bigg|\int_{\kappa(\Delta,s)}^{s}b(X^{\Delta}(\kappa(\Delta,r)),X^{\Delta}(\kappa(\Delta,r)-\tau),r)dr\bigg|ds\\
&+\E\int_{0}^{T\wedge\tau^{\Delta}(R)}\bigg|\int_{\kappa(\Delta,s)}^{s}\sigma(X^{\Delta}(\kappa(\Delta,r)),X^{\Delta}(\kappa(\Delta,r)-\tau),r)dB(r)\bigg|ds.
\end{split}
\end{equation}
By using  {\bf(C4)}, and the Burkholder-Davis-Gundy inequality, we can write that
\begin{equation}
\begin{split}
& \E\int_{0}^{T\wedge\tau^{\Delta}(R)}|p^{\Delta}(s)|ds \\
&\leq  \frac{1}{1-\kappa}\E\int_{0}^{T\wedge\tau^{\Delta}(R)}\int_{\kappa(\Delta,s)}^{s}|b(X^{\Delta}(\kappa(\Delta)),X^{\Delta}(\kappa(\Delta,r)-\tau),r)|drds\\
&+\frac{1}{1-\kappa}\E\int_{0}^{T\wedge\tau^{\Delta}(R)}4\sqrt{2}\E\bigg(\int_{\kappa(\Delta,s)}^{s}\|\sigma(X^{\Delta}(\kappa(\Delta)),X^{\Delta}(\kappa(\Delta,r)-\tau)r,)|^2dr\|\bigg)^{1/2}ds.
\end{split}
\end{equation}
Recalling the standing hypothesis {\bf(H)}, and then letting $\Delta\rightarrow0,$ we obtain that for all $T\in[0,\infty)$
\begin{equation*}
\E\int_{0}^{T\wedge\tau^{\Delta}(R)}|p^{\Delta}(s)|ds \to 0.
\end{equation*}
Therefore, the assumption {\bf(ii)} in the  localization lemma holds. Therefore for any $t\in[0,T]$ the localization lemma yields
$$
\sup_{0\le t \le T}|X^{\Delta_1}(t)-X^{\Delta_2}(t) | \to 0 \mbox{ in probability as } \Delta_1, \Delta_2 \to 0.
$$
By Lemma \ref{lemma3},  we then  have
\begin{equation}\label{UniformConvergence}
\sup_{t\in[0,T]}|X^{\Delta}(t)-X(t)|\xrightarrow{\mathbb{P}} 0 \quad  \mbox{as} \quad  \Delta\rightarrow 0.
\end{equation}
 To procede, we need to fix $T\in [0,\infty).$ By $(\ref{UniformConvergence})$ and the continuity of the path, we only need to show that the right hand side of $(\ref{eq23})$ converges almost surely to
\begin{equation*}
D(X(t-\tau))+X(0)-D(\xi(-\tau))+\int_{0}^{t}b(X(s),X(s-\tau),s)ds+\int_{0}^{t}\sigma(X(s),X(s-\tau),s)dB(s)
\end{equation*}
Since the uniform convergence is given in the equation $(\ref{UniformConvergence})$ on $[0,T],$ we also have
\begin{equation*}
\sup_{t\in [0,T]}|X^{\Delta}(\kappa(\Delta, t))-X(t)| \xrightarrow{\mathbb{P}} 0 \quad  \mbox{as} \quad  \Delta\rightarrow 0.
\end{equation*}
Let $Y^{\Delta}(t)= X^{\Delta}(\kappa(\Delta,t),t)$ and there exists a subsequence $(\Delta_{m})_{m \in \mathbb{N}}$ such that
\begin{equation*}
\sup_{t \in [0,T]}|Y^{\Delta_{m}}(t)-X(t)| \to 0 \quad  \quad a.s. \mbox{ when } m\rightarrow \infty. 
\end{equation*}
Moreover for $\bar Y(t):= \sup_{m\in \mathbb{N}} |Y^{\Delta_{m}}(t)|,$ we have
$
\sup_{t\in [0, T]}\bar Y(t)  < \infty \quad \mbox{a.s.}.
$
For the neutral term, it is easy to verify
\begin{equation}\label{neutral}
\begin{split}
|D(X(t-\tau))-D(Y^{\Delta_{m}}(\kappa(\Delta_{m},t)-\tau))|\leq k |X(t-\tau)-Y^{\Delta_{m}}(\kappa(\Delta_{m},t)-\tau)| \to 0 \quad \mbox{a.s.},
\end{split}
\end{equation}
as $m\rightarrow \infty.$
Define the $(\mathcal{F}_{t})$-stopping time
\begin{equation*}
\rho(R):=\inf\{t \in [0,T]:  \bar Y(t)>R  \} \wedge T
\end{equation*}
By   \eqref{bound}  and Lebesgue's dominated convergence theorem, on $\{t \leq \rho_{N}(R)\}$  we derive
\begin{equation}\label{eq25}
\int_0^tb(Y^{\Delta_{m}}(\kappa(\Delta_{m},s)),Y^{\Delta_{m}}(\kappa(\Delta_{m},s)-\tau),s)ds  \to \int_0^tb(X(s),X(s-\tau),s)ds \quad \mbox{a.s.}
\end{equation}
and
\begin{equation}\label{eq26a}
\lim_{m \to \infty}\int_{0}^{t}\E\left(\|\sigma(Y^{\Delta_{m}}(\kappa(\Delta_{m},s)),Y^{\Delta_{m}}(\kappa(\Delta_{m},s)-\tau),s)-\sigma(X(s),X(s-\tau),s)\|^{2}ds \right) = 0. 
\end{equation}
This  implies
\begin{equation}\label{eq26}
\int_{0}^{t}\sigma(Y^{\Delta_{m}}(Y^{\Delta_{m}}(\kappa(\Delta_{m},s)),Y^{\Delta_{m}}(\kappa(\Delta_{m},s)-\tau),s)dB(s)\to  \int_{0}^{t}\sigma(X(s),X(s-\tau))dB(s) \quad \mbox{a.s}.
\end{equation}
%$ The equation $(\ref{eq25})$
%shows that there exists a $\omega \in \Omega$,$R(\omega) \in [0,\infty]$
%such that  $\tau(R) = T$ for all $R \geq R(\omega)$ almost surely. \\
Due to the hypothesis {\bf(H)} for every $\omega \in \Omega$ there exists $N(\omega)\in [0,\infty)$ such that $\rho_{N}(R) = \rho(R)$ for all $N \geq N(\omega),$ so that
\begin{equation*}
\cup_{N\in\mathbb{N}}\{t\leq\rho_{N}(R) \} = \{t\leq \rho(R)\}.
\end{equation*}
This implies $(\ref{eq26})$ holds on $\{t\leq \rho(R)\}.$  However due to $
\sup_{t\in [0, T]}\bar Y(t)  < \infty \quad \mbox{a.s.}
$ for $\omega \in \Omega,$ there exists $R(\omega)\in [0,\infty)$ such that $\rho(R)=T$ for all $R \geq R(\omega).$ Hence, we have shown that all \eqref{neutral}, $(\ref{eq25})$ and $(\ref{eq26})$ hold almost surely. This completes the proof of existence.

For the uniqueness part, we suppose that $X(t)$ and $\bar{X}(t)$ are two solutions to \eqref{NSDDE} with the same initial data \eqref{initial}.  It is easy to see that the Euler numerical solution will converge to $X(t)$ and $\bar{X}(t)$, we must have
$$
\mathbb{P}\{\omega: X(t, \o)=\bar X(t, \o), 0\le t\le T\}=1.
$$
Therefore, we have complete the proof of uniqueness.

In order to estimate  the $p-$th moment, let $\hat{\kappa}$ be a negative number, we  define
\begin{equation*}
\rho(t) := \exp(\hat{\kappa}\alpha_{1}(t)) 
\end{equation*} 
By the It\^o formula and the assumption {\bf(C2)}, for all $t\in[0,T],$\\
\begin{equation}\label{Boundedness}
\begin{split}
&|X(t)-D(X(t-\tau))|^{2}\rho(t) \\
&=|\xi(0)-D(\xi(-\tau))|^2+\int_{0}^{t}\rho(s)\big[2\langle X(s)-D(X(s-\tau)),b(X(s),X(s-\tau),s)\rangle\\
&+||\sigma(X(s),X(s-\tau),s)||^2+\hat{\kappa}K_{1}(s)|X(s)-D(X(s-\tau))|^2\big]ds+M(t)\\
&\leq |\xi(0)-D(\xi(-\tau))|^2+\int_{0}^{t}\rho(s)K_{1}(s)\big[\hat{\kappa}|X(s)-D(X(s-\tau))|^2\\
&+K_{1}(s)(1+|X(s)|^{2})+\tilde{K}_{1}(s-\tau)(1+|X(s-\tau)|^{2})\big]ds+M(t), 
\end{split}
\end{equation} 
where  
\begin{equation*}
M(t) = \int_{0}^{t}\rho(s)2\langle X(s)-D(X(s-\tau)),\sigma(X(s),X(s-\tau),s)\rangle dB(s),
\end{equation*}
for all $t\in [0,T],$ which is a continuous local martingale with $M(0)=0.$ Noting that
$$
|X(s)-D(X(s-\tau))|\ge (1-\kappa )|X(s)|^2+\kappa (\kappa-1)|X(s-\tau)|^2,
$$
and $\r(s)$ is a non-increasing function, we have
\begin{equation}
\begin{split}
&|X(t)-D(X(t-\tau))|^{2}\rho(t) \\
&\leq |\xi(0)-D(\xi(-\tau))|^2+\int_{0}^{t}\rho(s)\big[K_{1}(s)\hat \kappa((1-\kappa)|X(s)|^{2}+\kappa(\kappa-1)|X(s-\tau)|^{2})\\
&+K_{1}(s)(1+|X(s)|^{2})+\tilde{K}_{1}(s-\tau)(1+|X(s-\tau)|^{2})\big]ds +M(t)\\
&\leq |\xi(0)-D(\xi(-\tau))|^2
+\int_{0}^{t}\rho(s)\big[K_{1}(s)(\bar{\kappa}((1-\kappa)+C_{1}(\tau)(\kappa^2-\kappa))+1)|X(s)|^{2}\big ]ds\\
&+\int_{0}^{t}\rho(s)\tilde{K}_{1}(s)|X(s)|^{2}ds+C_1(\tau)\int_{-\tau}^0\rho(\theta)K_{1}(\theta)\big[ \hat \kappa (\kappa^2-\kappa)\big]|\xi(\theta)|^2d \theta\\
&+\int_{-\tau}^0\rho(\theta)\tilde{K}_{1}(\theta)(1+|\xi(\theta)|^2)d \theta+\int_0^t\rho(s)(K_1(s)+\tilde{K}_{1}(s))ds +M(t).\\
\end{split}
\end{equation} 
Noting that $K_{1}(t) \geq \tilde{K}_{1}(t),$ then by choosing $\hat{\kappa}=-\frac{2}{(1-\kappa)+K_{1}(\tau)(\kappa^2-\kappa)},$  we can derive  
\begin{equation}\label{boundedness1}
\begin{split}
& \E|X(t)-D(X(t-\tau))|^{2} \rho(t)\leq \E|\xi(0)-D(\xi(-\tau))|^2\\
&+\E\int_{-\tau}^0\rho(\theta)K_{1}(\theta)\big[ \hat \kappa (\kappa^2-\kappa)\big]|\xi(\theta)|^2d \theta\\
&+\E\int_{-\tau}^0\rho(\theta)K_{1}(\theta)(1+|\xi(\theta)|^2)d \theta+\int_0^t2\rho(s)K_1(s)ds.\\
\end{split}
\end{equation}
An application of Lemma \ref{lemma2} yields that
\begin{equation*}
\begin{split}
&|X(t)|^{2} = |X(t)-D(X(t-\tau))+D(X(t-\tau))|^{2}\\
&\leq (1+\epsilon)(|D(X(t-\tau))|^{2}+\frac{|X(t)-D(X(t-\tau))|^{2}}{\epsilon})\\
&\leq (1+\epsilon)\kappa^{2}|X(t-\tau)|^{2}+\frac{1+\epsilon}{\epsilon}|X(t)-D(X(t-\tau))|^{2},
\end{split}
\end{equation*}
where assumption  {\bf(C4)} is applied.
For any $\kappa\in(0,1),$ let $\epsilon< \frac{1-\kappa^{2}}{\kappa^{2}},$ then take the expectation of both sides, finally take the supremum of both sides, we have
\begin{equation}\label{boundedness2}
\begin{split}
\sup_{0\leq t\leq T}\E|X(t)|^{2}\rho(t)\leq C\E||\xi||_{(-\tau,0)}^{2}+C\sup_{0\leq t\leq  T}\E|X(t)-D(X(t-\tau))|^{2}\rho(t).\\
\end{split}
\end{equation}
Since $\rho(t)$ is a positive function, and it is bounded for any $t\in[0,T],$ the required boundedness result follows by combining \eqref{boundedness1} and \eqref{boundedness2}. $\Box$

\section{Example}

In this section, we shall apply the Theorem \ref{EUtheorem} to the following nonlinear equation.
\begin{expl} {\rm
Consider an one-dimensional NSDDE, for any $k\in(-1,1),$ $t\in[0,T],$
\begin{equation}\label{example}
\begin{split}
d[X(t)-kX(t-\tau)] &= e^{c_{1}t}[1+X(t)-kX(t-\tau)-X^{3}(t)-k^{2}X(t)X^{2}(t-\tau)\\
&+kX^{2}(t)X(t-\tau)+k^{3}X^3(t-\tau)]ds+e^{c_{2}t}(1+X(t)-kX(t-\tau))dB(t),
\end{split}
\end{equation}
with the initial data $\{\xi(t):-\tau\leq t \leq 0\}$ $\in$ $C ([-\tau, 0];\mathbb{R}),$ where $B(t)$ is an one-dimensional Brownian motion and $c_1,c_2$ are two non-positive numbers with $c_1\leq c_2$.

It is not difficult  to verify that  {\bf(C1)}, {\bf(C4)} and {\bf(H)} hold. Now, we shall verify that  {\bf(C2)} and {\bf(C3)} hold respectively.
%Firstly, we note that $e^{ct}$ is locally integrable, and $e^{ct} = e^{c\tau}e^{c(t-\tau)},$ where $e^{c\tau}$ is a constant.\\
 We now compute
\begin{equation*}
\begin{split}
&e^{c_{1}t}\langle x-ky,1+x-ky-x^{3}-k^{2}xy^2+kx^{2}y+k^{3}y^{3}\rangle+e^{c_{2}t}(1+|x-ky|)^{2}\\
&=e^{c_{1}t}[(x-ky)+(x-ky)^{2}-(x-ky)^{2}(x^{2}+k^2y^{2})]+2e^{c_{2}t}|x-ky|^{2}+2e^{c_{2}t}\\
&\leq 4(e^{c_{1}t}+e^{c_{2}(t)})(1+|x|^{2})+ 4k^{2}(e^{c_{1}t}+e^{c_{2}t})(1+|y|^{2})\\
&\leq 4(e^{c_{1}t}+e^{c_{2}t})(1+|x|^{2})+ 4k^{2}(e^{c_{1}(t-\tau)}+e^{c_{2}(t-\tau)})(1+|y|^{2}).
\end{split}
\end{equation*}
Since $|k|<1,$ we have $4k^{2}(e^{c_{1}t}+e^{c_{2}t})< 4(e^{c_{1}t}+e^{c_{2}t})(1+|x|^{2}).$ This means {\bf(C2)} holds.
Moreover for all $|x|\vee |y|\vee |\bar{x}|\vee|\bar{y}| \leq R,$  
\begin{equation*}
\begin{split}
& e^{c_{1}t}\langle x-ky-\bar{x}+k\bar{y},x-ky-x^3-k^2xy^2+kx^{2}y+k^3y^{3}\\
&-\bar{x}+k\bar{y}+\bar{x}^{3}+k^{2}\bar{x}\bar{y}^{2}-k\bar{x}^{2}\bar{y}-k^3\bar{y}^3\rangle+e^{c_{2}t}|x-ky-\bar{x}+k\bar{y}|^{2}\\
&\leq e^{c_{1}t}\langle x-\bar{x}-k(y-\bar{y}),(x-\bar{x})-k(y-\bar{y})-(x-\bar{x})(x^{2}+x\bar{x}+\bar{x}^{2})\\
&-k^{2}(xy^{2}-x\bar{y}^{2}+x\bar{y}^{2}-\bar{x}\bar{y}^{2})+k(x^{2}y-x^{2}\bar{y}+x^{2}\bar{y}-\bar{x}^{2}\bar{y})\\
&+k^{3}(y-\bar{y})(y^2+y\bar{y}+\bar{y}^{2})\rangle+e^{c_{2}t}|x-ky-\bar{x}+k\bar{y}|^{2}\\
&\leq \langle x-\bar{x}-k(y-\bar{y}),(x-\bar{x})-k(y-\bar{y})-(x-\bar{x})(x^{2}+x\bar{x}+\bar{x}^{2})\\
&-k^{2}(x(y-\bar{y})(y+\bar{y})+(x-\bar{x})\bar{y}^{2})+k(x^{2}(y-\bar{y})+(x-\bar{x})(x+\bar{x})\bar{y})\\
&+k^{3}(y-\bar{y})(y^2+y\bar{y}+\bar{y}^{2})\rangle+2e^{c_{2}t}(|x-\bar{x}|^2+k^{2}|y-\bar{y}|^2)\\
&\leq (1+\bar{y}^{2}+(x+\bar{x})\bar{y}+P(x,y,\bar{x},\bar{y}))e^{c_{1}t}|x-\bar{x}|^{2}\\
&+(k^{2}+k^3x(y+\bar{y})+P(x,y,\bar{x},\bar{y}))e^{c_{1}t}|y-\bar{y}|^{2}\\
&+2e^{c_{2}t}(|x-\bar{x}|^2+k^{2}|y-\bar{y}|^2),
\end{split}
\end{equation*}
where 
\begin{equation*}
\begin{split}
& P(x,y,\bar{x},\bar{y}) = \bigg(|2k|+|k^{2}x(y+\bar{y})|+|kx^2|+|k^{3}(y^2+y\bar{y}+\bar{y}^{2})|\\
&+|k(x^{2}+x\bar{x}+\bar{x}^{2})|+|k^{3}\bar{y}^{2}|+|k^{2}(x+\bar{x})\bar{y}|\bigg)^{1/2}.\\
\end{split}
\end{equation*}
Noting that $|x|\vee |y|\vee |\bar{x}|\vee|\bar{y}| \leq R,$ we obtain
\begin{equation*}
\begin{split}
& e^{c_{1}t}\langle x-ky-\bar{x}+k\bar{y},x-ky-x^3-k^2xy^2+kx^{2}y+k^3y^{3}\\
&\le  (1+R^{2}+(R+R)R+P(R,R,R,R))|x-\bar{x}|^{2}\\
&+(k^{2}+k^3R(R+R)+P(R,R,R,R))|y-\bar{y}|^{2}
+2(|x-\bar{x}|^2+k^{2}|y-\bar{y}|^2)\\
&=: C(R)|x-\bar{x}|^{2}+\tilde C(R)|y-\bar{y}|^2),
\end{split}
\end{equation*}
  we can see $C(R)\ge \tilde C(R).$ This means {\bf(C3)} holds.
  
  Hence, by the Theorem \ref{EUtheorem} the equation \eqref{example} has a unique solution.}
\end{expl}

\end{document}